\documentclass[11pt]{article}
\usepackage{amsmath,amsthm,amsfonts}

\newcommand{\LL}{\mathcal L}
\newcommand{\M}{\mathcal M}
\renewcommand{\=}{\doteq}

\newtheorem{thm}{Theorem}[section]
 \newtheorem{cor}[thm]{Corollary}
 \newtheorem{prop}[thm]{Proposition}
 \newtheorem{lemma}[thm]{Lemma}
\theoremstyle{definition}
 \newtheorem{defn}[thm]{Definition}
\theoremstyle{definition}
 \newtheorem{rem}[thm]{Remark}
\numberwithin{equation}{section}
\begin{document} 
\title{\bf Moufang symmetry II.\\ Moufang-Mal'tsev pairs and triality}
\author{Eugen Paal}
\date{}
\maketitle

\thispagestyle{empty}

\begin{abstract}
A concept of the Moufang-Malt'tsev pair is elaborated. This concept is based on the generalized Maurer-Cartan equations of  a local analytic Moufang loop. Triality can be seen as a fundamental property of such pairs. Based on  triality, the Yamagutian is constructed. Properties of the Yamagutian are studied.
\end{abstract}

\section{Introduction}

In \cite{Paal}, the generalized Maurer-Cartan equations for  a local analytic Moufang loop were found.   In this paper, we elaborate a concept of the \emph{Moufang-Mal'tsev pair} which is based on the generalized Maurer-Cartan equations. Such a concept may also be inferred from the theory of  \emph{alternative} algebras \cite{Schafer52,Schafer66}. Triality can be seen as a fundamental property of such pairs. Based on triality, the Yamagutian is constructed. Properties of the Yamagutian are studied.

\section{Moufang-Mal'tsev pairs}

Let $\M$ be an anti-commutative algebra and let $\LL$ be a Lie algebra. Throughout this paper we assume that both  algebras  have the same base field ${\mathbb F}$ ($\= {\mathbb R}$ or ${\mathbb C}$;  as a matter of fact, only $\operatorname{char}{\mathbb F}\ne 2, 3$ is essential).  Denote by $(S,T)$ a pair of the linear maps 
$S,T: \M\to \LL$. 

\begin{defn}[Moufang-Mal'tsev pair]
We call $(S,T)$ a  \emph{Moufang-Mal'tsev pair} if
\begin{subequations}
\label{m-m}
\begin{align} 
[S_{x},S_{y}]&=S_{[x,y]}-2[S_{x},T_{y}]\\
[T_{x},T_{y}]&=T_{[y,x]}-2[T_{x},S_{y}]
\end{align}
\end{subequations}
for all $x,y$ in $\M$.  
\end{defn}
We call (\ref{m-m}a,b) the Moufang-Mal'tsev relations. Note that the same brackets $[\cdot,\cdot]$ are used to denote multiplications in $\M$ and $\LL$.

\begin{prop} 
Let $(S,T)$ a Moufang-Mal'tsev pair. Then
\[
[S_{x},T_{y}]=[T_{x},S_{y}]
\] 
for all $x,y$ in $\M$.
\end{prop}

\begin{proof} 
Use anti-commutativity in $\LL$ and $\M$.
\end{proof}

\begin{cor}[minimality conditions]
The Moufang-Mal'tsev conditions read
\[
2[S_{x},T_{y}]
=S_{[x,y]}-[S_{x},S_{y}]
=T_{[y,x]}-[T_{x},T_{y}]
=2[T_{x},S_{y}]
\]
\end{cor}

\begin{rem}
Deviations of $S$ and $T$ from algebra  maps may be considered as "minimal" for the Moufang-Mal'tsev pair.
\end{rem}

\section{Triality}

\begin{prop}
\label{-s-t}
Let $(S,T)$ be a Moufang-Mal'tsev pair. 
Then $(-T,-S)$ is a Moufang-Mal'tsev pair as well.
\end{prop}

\begin{proof} 
By using anti-commutativity in $M$, rewrite the Moufang-Mal'tsev conditions of $(S,T)$ as  follows:
\begin{align*}
[-T_{x},-T_{y}]&=-T_{[x,y]}-2[-T_{x},-S_{y}]\\
[-S_{x},-S_{y}]&=-S_{[y,x]}-2[-S_{x},-T_{y}]\tag*{\qed}
\end{align*}
\renewcommand{\qed}{}
\end{proof}

Consider a triple $(S,T,P)$ of the linear maps $S,T,P:M\to L$, such that 
\begin{equation*}
S+T+P=0
\end{equation*}
With a given triple $(S,T,P)$ we can associate the pairs 
\begin{equation}
\label{stp}
(S,T),\quad (-T,-S),\quad (T,P),\quad (-P,-T),\quad (P,S),\quad (-S,-P)
\end{equation}  

\begin{lemma}[triality]  
\label{triality}
Let some pair from (\ref{stp}) be a  Moufang-Mal'tsev pair. Then all other pairs from (\ref{stp}) are  the Moufang-Mal'tsev pairs as well. 
\end{lemma}

\begin{proof}
Assume that $(S,T)$ is a Moufang-Mal'tsev pair. Then we know from above Proposition \ref{-s-t} that $(-T,-S)$ is a Moufang-Mal'tsev pair. The required minimality conditions for $(T,P)$ and $(P,S)$ read, respectively,
\begin{align*}
2[T_{x},P_{y}]
&=T_{[x,y]}-[T_{x},T_{y}]
=P_{[y,x]}-[P_{x},P_{y}]
=2[P_{x},T_{y}]\\
2[P_{x},S_{y}]
&=P_{[x,y]}-[P_{x},P_{y}]
=S_{[y,x]}-[S_{x},S_{y}]
=2[S_{x},P_{y}]
\end{align*}
As an example, calculate
\begin{align*}
-2[S_{x},P_{y}]&=2[P_{x}+T_{x},P_{y}]\\
    &=2[P_{x},P_{y}]+2[P_{x},T_{y}]\\
    &=[P_{x},P_{y}]+[P_{x},P_{y}]+2[P_{x},T_{y}]\\
    &=[P_{x},P_{y}]+[S_{x}+T_{x},S_{y}+T_{y}]-2[S_{x}+T_{x},T_{y}]\\
    &=[P_{x},P_{y}]+[S_{x},S_{y}]+[S_{x},T_{y}]+[T_{x},S_{y}]
        +[T_{x},T_{y}]-2[S_{x},T_{y}]-2[T_{x},T_{y}]\\
    &=[P_{x},P_{y}]+[S_{x},S_{y}]-[T_{x},T_{y}]\\
    &=[P_{x},P_{y}]+S_{[x,y]}-2[S_{x},T_{y}]
        -T_{[y,x]}+2[T_{x},S_{y}]\\
    &=[P_{x},P_{y}]-P_{[x,y]}.
\end{align*}
All other required equalities in (\ref{stp}) can be verified  in the same way. The result tells us that $(T,P)$ and   
$(P,S)$ are the Moufang-Mal'tsev pairs, which in turn implies that $(-P,-T)$ and $(-S,-P)$ must be the Moufang-Mal'tsev pairs as well.
\end{proof}

It follows from Lemma \ref{triality} that the Moufang-Mal'tsev pair $(S,T)$ is invariant under the  substitutions
\begin{align*}
\sigma   &\=  (S\to  T\to P\to S)\\
\tau     &\=  (S\to -T\to S)(P\to -P)\\
\sigma^{2}          &=      (S\to  P\to T\to S)\\
\sigma\circ\tau     &=      (S\to -P\to S)(T\to -T)\\
\sigma^{2}\circ\tau &=      (T\to -P\to T)(S\to -S)
\end{align*}
which we call the \emph{triality substitutions}.  So it is natural to proclaim:

\begin{thm}[principle of triality]
All algebraic consequences of the Moufang-Mal'tsev conditions are triality invariant. 
\end{thm}

Such a symmetry we call \emph{triality}. It suggests that we should try to handle the Moufang-Mal'tsev 
pairs in triality symmetric manner. 

In particular, by using the \emph{triality conjugation},
\begin{align*}
 P^{+}&\= S-T=P+2S=-P-2T\\
 S^{+}&\=  T-P=S+2T=-S-2P\\
 T^{+}&\=  P-S=T+2P=-T-2S
\end{align*}
with the evident property 
\begin{equation*} 
S^{+}+T^{+}+P^{+}=0
\end{equation*} 
one can rewrite the Moufang-Mal'tsev conditions as follows: 
\begin{align*}   
[S_{x},S^{+}_{y}]&=[S^{+}_{x},S_{y}]=S_{[x,y]}\\
[T_{x},T^{+}_{y}]&=[T^{+}_{x},T_{y}]=T_{[x,y]}\\
[P_{x},P^{+}_{y}]&=[P^{+}_{x},P_{y}]=P_{[x,y]}
\end{align*}
Note that  
\begin{align*}
3P&=T^{+}-S^{+}=P^{+}+2T^{+}=-P^{+}-2S^{+}\\ 
3S&=P^{+}-T^{+}=S^{+}+2P^{+}=-S^{+}-2T^{+}\\
3T&=S^{+}-P^{+}=T^{+}+2S^{+}=-T^{+}-2P^{+}
\end{align*}
which means that the triality conjugation is invertible.

\section{Yamagutian}

We introduced the triple $(S,T,P)$ via the triality symmetric identity (\ref{stp}). Following triality, it is natural to search for other but nontrivial triality invariant combinations of the maps from the triple $(S,T,P)$. \begin{defn}[Yamagutian \cite{Yam}] 
The \emph{Yamagutian} of    $(S,T)$ is the skew-symmetric bilinear map $Y: M\otimes M\to L$ defined (cf (\ref{stp})) by
\begin{align*}
6Y(x;y)
&\=  [S_{x},S_{y}]+[T_{x},T_{y}]+[P_{x},P_{y}]\\
&=-Y(y;x)
\end{align*}
\end{defn}
We can see the evident but important
\begin{prop}
The Yamagutian  $Y$ is triality invariant. 
\end{prop}

By triality symmetry, the Yamagutian $Y$ can be redefined in several useful ways. In particular,
\begin{align*}
6Y(x;y) 
&=3[S_{x},S_{y}]-S^{+}_{[x,y]}\\
&=3[T_{x},T_{y}]-T^{+}_{[x,y]}\\
&=3[P_{x},P_{y}]-P^{+}_{[x,y]}\\
\intertext{and one can also verify that}
6Y(x;y) 
&=2P^{+}_{[x,y]}-6[S_{x},T_{y}]\\
&=2S^{+}_{[x,y]}-6[T_{x},P_{y}]\\
&=2T^{+}_{[x,y]}-6[P_{x},S_{y}]
\end{align*}
Later we shall need the 
\begin{prop}
Let $(S,T)$ be a Moufang-Mal'tsev pair. Then
\begin{subequations}
\label{stp+}
\begin{align}
6Y(x;y)
&=[S^{+}_{x},S^{+}_{y}]+S^{+}_{[x,y]}\\
&=[T^{+}_{x},T^{+}_{y}]+T^{+}_{[x,y]}\\
&=[P^{+}_{x},P^{+}_{y}]+P^{+}_{[x,y]}
\end{align}
\end{subequations}
for all $x,y$ in $M$.
\end{prop}

\begin{proof} 
Due to triality, check only the last formula (\ref{stp+}c):
\begin{align*}
[P^{+}_{x},P^{+}_{y}]
&=[S_{x}-T_{x},S_{y}-T_{y}]\\
&=[S_{x},S_{y}]-2[S_{x},T_{y}]+[T_{x},T_{y}]\\
&= \frac{1}{3}Y(x;y)+\frac{1}{3}S^{+}_{[x,y]}
+\frac{1}{3}Y(x;y)-\frac{2}{3}P^{+}_{[x,y]}
+\frac{1}{3}Y(x;y)+\frac{1}{3}T^{+}_{[x,y]}\\
&=Y(x;y)-P^{+}_{[x,y]}\tag*{\qed}
\end{align*}
\renewcommand{\qed}{}
\end{proof}

\begin{rem}
Formulae (\ref{stp+}a--c) tell us that the Yamagutian $Y$ measures the deviation of $S^{+}$, $T^{+}$ and $P^{+}$ from the anti-algebra maps. 
\end{rem}

\begin{cor}
We have
\begin{equation*}
18Y(x,y)= [S^+_{x},S^+_{y}]+[T^+_{x},T^+_{y}]+[P^+_{x},P^+_{y}]
\end{equation*}
\end{cor}

\begin{thm}
Let $(S,T)$ be a Moufang-Mal'tsev pair. Then
\begin{align*}
[S_{x},S_{y}]&=2Y(x;y)+\frac{1}{3}S_{[x,y]}+\frac{2}{3}T_{[x,y]}\\
[S_{x},T_{y}]&=-Y(x;y)+\frac{1}{3}S_{[x,y]}-\frac{1}{3}T_{[x,y]}\\
[T_{x},T_{y}]&=2Y(x;y)-\frac{2}{3}S_{[x,y]}-\frac{1}{3}T_{[x,y]}
\end{align*}
\end{thm}

\begin{proof}
Evident.
\end{proof}

\begin{cor}
By triality, we have
\begin{align*}
[T_{x},T_{y}]&=2Y(x;y)+\frac{1}{3}T_{[x,y]}+\frac{2}{3}P_{[x,y]}\\
[T_{x},P_{y}]&=-Y(x;y)+\frac{1}{3}T_{[x,y]}-\frac{1}{3}P_{[x,y]}\\
[P_{x},P_{y}]&=2Y(x;y)-\frac{2}{3}T_{[x,y]}-\frac{1}{3}P_{[x,y]}
\end{align*}
and
\begin{align*}
[P_{x},P_{y}]&=2Y(x;y)+\frac{1}{3}P_{[x,y]}+\frac{2}{3}S_{[x,y]}\\
[P_{x},S_{y}]&=-Y(x;y)+\frac{1}{3}P_{[x,y]}-\frac{1}{3}S_{[x,y]}\\
[S_{x},S_{y}]&=2Y(x;y)-\frac{2}{3}P_{[x,y]}-\frac{1}{3}S_{[x,y]}
\end{align*}
\end{cor}

\begin{prop} 
Let $(S,T)$ be a Moufang-Mal'tsev pair. Then 
\begin{subequations}
\label{y-stp}
\begin{align}
6[Y(x;y),S_{z}]&=3[[S_{x},S_{y}],S_{z}]-S_{[[x,y],z]}\\
6[Y(x;y),T_{z}]&=3[[T_{x},T_{y}],T_{z}]-T_{[[x,y],z]}\\
6[Y(x;y),P_{z}]&=3[[P_{x},P_{y}],P_{z}]-P_{[[x,y],z]}
\end{align}
\end{subequations}
for all $x,y,z$ in $M$.
\end{prop}
\begin{proof}  
Due to triality, only the first identity must be checked:
\begin{align*}
6[Y(x;y),S_{z}]
&=[3[S_{x},S_{y}]-S^{+}_{[x,y]},S_{z}]\\
&=3[[S_{x},S_{y}],S_{z}]-[S^{+}_{[x,y]},S_{z}]\\
&=3[[S_{x},S_{y}],S_{z}]-S_{[[x,y],z]}\tag*{\qed}
\end{align*}
\renewcommand{\qed}{}
\end{proof}

\begin{cor}
Adding formulae (\ref{y-stp}a--c) we obtain (cf (\ref{stp})) the triality symmetric identity  
\begin{equation*}
[[S_{x},S_{y}],S_{z}]+[[T_{x},T_{y}],T_{z}]+[[P_{x},P_{y}],P_{z}]=0
\end{equation*}
\end{cor}

\begin{cor}
In (\ref{y-stp}a--c) make twicely cyclic permutation of $x,y,z$  and add the resulting equalities with the original ones. Then we obtain
\begin{align*}
6[Y(x;y),S_{z}]+6[Y(y,z),S_{x}]+6[Y(z;x),S_{y}]&=S_{J(x,y,z)}\\
6[Y(x;y),T_{z}]+6[Y(y,z),T_{x}]+6[Y(z;x),T_{y}]&=T_{J(x,y,z)}\\
6[Y(x;y),P_{z}]+6[Y(y,z),P_{x}]+6[Y(z;x),P_{y}]&=P_{J(x,y,z)}
\end{align*}
\end{cor}

\begin{prop}
Let $(S,T)$ be a Moufang-Mal'tsev pair. Then
\begin{subequations}
\label{stp++}
\begin{align}
[S^{+}_{[x,y]},S^{+}_{z}]+[S^{+}_{[y,z]},S^{+}_{x}] +[S^{+}_{[z,x]},S^{+}_{y}]
&=P^{+}_{J(x,y,z)}\\
[T^{+}_{[x,y]},T^{+}_{z}]+[T^{+}_{[y,z]},T^{+}_{x}] +[T^{+}_{[z,x]},T^{+}_{y}]
&=T^{+}_{J(x,y,z)}\\
[P^{+}_{[x,y]},P^{+}_{z}]+[P^{+}_{[y,z]},P^{+}_{x}] +[P^{+}_{[z,x]},P^{+}_{y}]
&=P^{+}_{J(x,y,z)}
\end{align}
\end{subequations}
for all $x,y,z$ in $M$. 
\end{prop}

\begin{proof} 
Subtracting (\ref{y-stp}b) from (\ref{y-stp}a, we obtain
\begin{equation*}
6[Y(x,y),P^{+}_{z}]
=3[[S_{x},S_{y}],S_{z}]-3[[T_{x},T_{y}],T_{z}]-P^{+}_{[[x,y],z]}
\end{equation*}
On the other hand, using (\ref{y-stp}c), we have 
\begin{equation*}
6[Y(x;y),P^{+}_{z}]
=[[P^{+}_{x},P^{+}_{y}],P^{+}_{z}]+[P^{+}_{[x,y]},P^{+}_{z}]
\end{equation*}
and so we obtain
\begin{equation*}
[P^{+}_{[x,y]},P^{+}_{z}]
=-P^{+}_{[[x,y],z]}
+3[[S_{x},S_{y}],S_{z}]
-3[[T_{x},T_{y}],T_{z}]
-[[P^{+}_{x},P^{+}_{y}],P^{+}_{z}]
\end{equation*}
Now make twicely cyclic permutation of $x,y,z$ and add the resulting equalities with the original one. 
Then, using Jacobi conditions (in $L$) and the definition of $J$ on the right hand-side of the resulting equality we obtain (\ref{stp++}c). The remaining identities (\ref{stp++}a,b) are evident from triality.
\end{proof}

\begin{lemma} 
Let $(S,T)$ be a Moufang-Mal'tsev pair. Then 
\begin{equation}
\label{y-xyz}
Y([x,y];z)+Y([y,z];x)+Y([z,x];y)=0,\quad \forall x,y,z \in M
\end{equation}
\end{lemma}

\begin{proof}  
Use (\ref{stp+}c) to obtain   
\begin{equation*}
6Y([x,y];z)=[P^{+}_{[x,y]},P^{+}_{z}]+P^{+}_{[[x,y],z]}
\end{equation*}
Make here twicely the cyclic permutation of $x,y,z$ and add the resulting equalities with the original one. Then use (\ref{stp++}c) and the definition of $J$ to obtain the desired identity (\ref{y-xyz}).
\end{proof}

\section*{Acknowledgement}

Research was in part supported by the Estonian Science Foundation, Grant 6912.

\bigskip\noindent
Department of Mathematics\\
Tallinn University of Technology\\
Ehitajate tee  5, 19086 Tallinn, Estonia\\
E-mail: eugen.paal@ttu.ee

\end{document}